\documentclass{amsart}

\usepackage{amsthm,amssymb,latexsym,amsmath}
\usepackage{mathrsfs}
\usepackage{graphicx}
\usepackage{color}
\usepackage[brazil]{babel}
\usepackage[latin1]{inputenc}
\usepackage{graphicx}%
\usepackage{multirow}%
\usepackage{amsmath,amssymb,amsfonts}%
\usepackage{amsthm}%
\usepackage{physics}
\usepackage{mathrsfs}%
\usepackage[title]{appendix}%
\usepackage{xcolor}%
\usepackage{textcomp}%
\usepackage{manyfoot}%
\usepackage{booktabs}%
\usepackage{algorithm}%
\usepackage{algorithmicx}%
\usepackage{algpseudocode}%
\usepackage{listings}%
\usepackage{amsthm,amssymb,latexsym,amsmath}
\usepackage{mathrsfs}
\usepackage{graphicx}
\usepackage{color}
\usepackage[all]{xy}
\usepackage{graphicx}%
\usepackage{multirow}%
\usepackage{amsmath,amssymb,amsfonts}%
\usepackage{amsthm}%
\usepackage{mathrsfs}%
\usepackage[title]{appendix}%
\usepackage{xcolor}%
\usepackage{textcomp}%
\usepackage{manyfoot}%
\usepackage{booktabs}%
\usepackage{algorithm}%
\usepackage{algorithmicx}%
\usepackage{algpseudocode}%
\usepackage{listings}%

\newcommand{\CC}{{\mathbb C}}

\newcommand{\OO}{\mathcal O}

\newcommand{\fol}{\mathcal{F}}

\newcommand{\NN}{{\mathbb{N}}}

\newcommand{\PP}{\mathbb{P}}

\newcommand{\Sing}{{\rm Sing}}

\newtheorem{lema}{Lemma}[section]
\newtheorem{cor}[lema]{Corollary}
\newtheorem{teo}[lema]{Theorem}
\newtheorem*{teo1234}{Theorem}
\newtheorem*{teo1}{Theorem 1}
\newtheorem*{teo2}{Theorem 2}

\newtheorem{prop}[lema]{Proposition}

\theoremstyle{definition}
\newtheorem{remark}[lema]{Remark}

\newtheorem{exe}[lema]{Example}

\begin{document}

\title{Formulas for Residues of Type Camacho-Sad and Applications}

\hyphenation{folia-tion}

\begin{abstract}
In this paper, we provide formulas for the sum of residues of type Camacho-Sad of a holomorphic foliation with respect to an invariant analytic subvariety. As application, in context of projective foliations, we obtain a formula that relates the sum these residues with the degree and other characteristics of the invariant subvariety. 
Furthermore, we establish sufficient conditions ensuring that an irreducible curve is invariant by a projective foliation on $\PP^2$. In addition, we provide an adjunction formula for hypersurfaces with non-isolated singularities and obtain explicit formulas for the Milnor number of these hypersurfaces.
 
\end{abstract}

\author{Diogo da Silva Machado}

\address{\noindent  Diogo da Silva Machado\\
Departamento de Matem\'atica \\
Universidade Federal de Vi\c cosa\\
Avenida Peter Henry Rolfs, s/n - Campus Universitário \\
36570-900 Vi\c cosa- MG, Brazil} \email{diogo.machado@ufv.br}




\maketitle



\section{Introduction}

The Separatrix Theorem, proved by C.~Camacho and P.~Sad in~\cite{cama}, is a classical result in foliation theory. It guarantees the existence of a local integral curve (a separatrix) for any holomorphic vector field with an isolated singularity at the origin of the complex plane. These authors defined the local index of a holomorphic vector field relative to a nonsingular invariant curve, which was the main tool used in their proof. This index is known in the literature as the \emph{Camacho-Sad index}, and many generalizations have since been developed.

In the case where the invariant curve is a singular curve on the projective plane $\PP^2$, a generalization was given by A.~Lins Neto~\cite{AL1}. When the ambient surface is a complete intersection in \( \mathbb{P}^n \), a generalization was obtained by M.~G.~Soares~\cite{ccc0}. T.~Suwa, in~\cite{Suw}, further extended the definition of the index to allow invariant curves with singularities. He proved that for a compact curve invariant by a one-dimensional singular foliation on a surface, the sum of the indices equals the self-intersection number of the curve. In higher dimensions, when the invariant subvariety is nonsingular, the index has been considered by B.~Gmira~\cite{gmi} and A.~Lins Neto~\cite{AL2} in the case of invariant hypersurfaces, and in the general case by D.~Lehmann~\cite{Lehm0}. 

Subsequently, in arbitrary dimension and for an invariant singular subvariety, D. Lehmann and T. Suwa \cite{LS00},\cite{Suw0}  obtained the index --- also referred to as the {\it residue of type Camacho-Sad}. By employing a suitable adaptation of Chern-Weil theory to \v{C}ech-de Rham cohomology, they extended the Baum-Bott residue theory~\cite{BB1} to the setting of foliations restricted to their invariant subvarieties. In this framework, the index is obtained as the localization of the normal bundle associated with the invariant subvariety.

More recently, F.~Sancho de Salas~\cite{sancho} obtained a generalization of the Camacho-Sad index via a purely algebraic approach, using Grothendieck-type residues associated with a Pfaff system relative to an invariant subscheme. His construction is valid over an arbitrary base scheme and in any characteristic.

It is important to note that the computation and determination of residues of this nature is, in general, a difficult task, which underscores the significance of the results presented in this article.

In the present work, we establish a residue formula (see Theorem 1) for computing the sum of residues of type Camacho-Sad (Camacho-Sad index) of a $p$-dimensional  foliation $ \fol $ relative to an invariant analytic subvariety $ V $, of arbitrary dimension. We consider the case where $V$ is a local complete intersection defined as the zero set of a section of a vector bundle (this is the case, for example, of hypersurfaces, global complete intersections and complete intersections in projective space).

As an application of Theorem 1, we consider foliations on $\PP^n$ that leave invariant a (possibly singular) hypersurface, and we prove a formula relating the sum of Camacho-Sad type residues to the multidegree and other characteristics of the invariant hypersurface (see Theorem 2). This theorem generalizes a result from \cite{AL1}.

Proposition \ref{teo0} provides an adjunction formula for hypersurfaces with non-isolated singularities. Theorem \ref{parv} and Corollary \ref{parvv} provide explicit formulas for the Milnor number of these hypersurfaces.

Also as an application, in the context of holomorphic foliations on $\PP^2$, Corollary \ref{corvv} and Corollary \ref{corrr} establish conditions under which an irreducible curve is invariant by a foliation. In particular, Corollary \ref{corrr} recovers a result due to D. Cerveau and A. Lins Neto \cite[Ch. 3, Theorem 1]{Alc}.

\section{Statement of Theorem 1}

Recall the definition of the residue of type Camacho-Sad \cite[Ch.VI, 6]{Suw2}. Let $X$ be an $n$-dimensional complex compact manifold and $\fol$ a $p$-dimensional singular holomorphic foliation on $X$ which leaves invariant a local complete intersection $V$, of codimension $k$, defined by a section of a holomorphic vector bundle $N$ on $X$. Let $S_{\lambda}$ be a compact connected component of $S(\fol,V):= (Sing(\fol)\cap V)\cup Sing(V)$ and $U$ an open neighborhood of $S_{\lambda}$ in $V$ such that $U_0:= U- S_{\lambda}$ is in $V - S(\fol,V)$. Also, let $\widetilde{U_1}$ be a regular neighborhood of $S_{\lambda}$ in $X$ with $\widetilde{U_1}\cap S(\fol,V) \subset U$ and $\widetilde{U_0}$ a tubular neighborhood of $U_0$ in $X$ with projection $\rho: \widetilde{U_0} \longrightarrow U_0$. By taking the covering $\mathcal{U} = \{\widetilde{U_0},\widetilde{U_1}\}$ of $\widetilde{U} = \widetilde{U_0} \cup \widetilde{U_1}$, we consider the associated $\check{\mbox{C}}$ech-de Rham complex $A^{\bullet}(\mathcal{U})$.
  
Let \( F \) be the vector bundle over \( X \setminus \operatorname{Sing}(\mathcal{F}) \) whose sheaf of germs of holomorphic sections, \( \mathcal{O}(F) \), coincides with \( \left.\mathcal{F}\right|_{X \setminus \operatorname{Sing}(\mathcal{F})} \). It follows from \cite[Ch.~VI, Lemma~6.6]{Suw2} that the restriction of the bundle \( N \) to \( V \setminus S(\mathcal{F}, V) \) is an \( F \)-bundle. Moreover, since the restriction \( \left.N\right|_{\widetilde{U_0}} \) is diffeomorphic to the pullback \( \rho^{\ast}N_{U_0} \) of the normal bundle \( N_{U_0} \) of the (manifold) \( U_0 \), we can regard the connections on \( N_{U_0} \) as \( \left.F\right|_{U_0} \)-connections.

Let $\nabla'_{0}$ be an \( \left.F\right|_{U_0} \)-connection for \( N_{U_0} \) and $\nabla'_1$ an arbitrary connection for $N$ on $\widetilde{U_1}$, then for a homogeneous symmetric polynomial $\varphi$ of degree $d$ the characteristic class $\varphi(N)$ is represented by the cocycle
$$\varphi(\nabla'_{\ast}) = (\varphi(\nabla'_{0}), \varphi(\nabla'_{1}), \varphi(\nabla'_{0}, \nabla'_{1} ))$$ 
\noindent in $\check{\mbox{C}}$ech-de Rham cohomology group
$H^{2n}(A^{\bullet}(\mathcal{U}))$. Since $\nabla'_{0}$ is an $\eval{F}_{U_0}$-connection, it follows from the Bott vanishing theorem \cite[Ch.II, Theorem 9.11]{Suw2} that
$$
\varphi(\nabla'_{0}) \equiv 0,
$$
\noindent if $d>n-k-p$. Consequently, the cocycle $\varphi(\nabla'_{\ast})$ is in the relative  $\check{\mbox{C}}$ech-de Rham complex $A^{2d}(\mathcal{U}, \widetilde{U_0})$ and it defines a class in the relative cohomology $H^{2d}(U, U - S_{\lambda}; \CC)$, denoted by $\varphi_{S_{\lambda}}(N,\fol)$. The residue (of type Camacho-Sad), denoted by $\mbox{Res}_{\varphi}(\fol,N,  S_{\lambda})$, is defined as the image of that class by Alexander duality $H^{2d}(U, U - S_{\lambda}; \CC)  \simeq H_{2(n-k)-2d}(S_{\lambda};\CC)$ (see \cite[Ch.VI, Theorem 6.9]{Suw2}).

\begin{teo1} Let $X$ be an $n$-dimensional complex compact manifold and $\fol$ a $p$-dimensional singular holomorphic foliation on $X$ which leaves invariant a local complete intersection $V$, of codimension $k$, defined by a section of a holomorphic vector bundle $N$ on $X$. The following hold: 
\begin{itemize}
\item [(a)] $\displaystyle
\sum_{\lambda=1}^r \mbox{Res}_{\varphi}(\fol,N,S_{\lambda}) = [X] \frown \left(\varphi(N) \wedge c_k(N)\right) \in H_{0}(M;\CC) \simeq \CC,$\\\\
\noindent where $S_1, \ldots, S_r$ are all compact connected components of $S(\fol,V)$, $\varphi$ is any homogeneous symmetric polynomial of degree $d=n-k$ and $H_{0}(M;\CC)$ denotes the $0$-th singular homology of $M$.
\medskip
\item [(b)] If $V = \bigcap_{l =1}^kD_{l}$ is a complete intersection  of $k$ hypersurfaces $D_{l}$, then
\\
{\small{\begin{itemize}
\item [(i)] $\displaystyle
\sum_{\lambda=1}^r \mbox{Res}_{c_{n-k}}(\fol,N,S_{\lambda}) = [X] \frown \hspace{-0.3cm} \sum_{1\leq l_1 <\ldots <l_{n-k} \leq k}\hspace{-0.8cm} c_1(L(D_{l_1}))... c_1(L(D_{l_{n-k}}))\prod^k_{l=1}c_1(L(D_l)),$
\medskip
\item [(ii)] $\displaystyle
\sum_{\lambda=1}^r \mbox{Res}_{c_1^{n-k}}(\fol,N,S_{\lambda}) = [X] \frown  \left(\sum_{l=1}^{k} c_l(L(D_{l}))\right)^{n-k}\prod^k_{l=1}c_1(L(D_l)),$ 
\end{itemize}
}}
\noindent where $L(D_{l})$ denotes the line bundle on $X$ associated to the divisor $D_{l}$.
\medskip
\item [(c)] In particular, if $X=\PP^n$ and $deg(D_{l}) = d_{l}$, for each $l=1,\ldots, k$, then
\\
\begin{itemize}
\item [(i)] $\displaystyle
\sum_{\lambda=1}^r \mbox{Res}_{c_{n-k}}(\fol,N,S_{\lambda}) = \hspace{-0.3cm} \sum_{1\leq l_1 <\ldots <l_{n-k} \leq k}\hspace{-0.8cm} (d_{l_1}\ldots d_{l_{n-k}})(d_{1}\ldots d_{k}),$
\medskip
\item [(ii)] $\displaystyle
\sum_{\lambda=1}^r \mbox{Res}_{c_1^{n-k}}(\fol,N,S_{\lambda}) = \left(d_1 +\ldots + d_k\right)^{n-k}(d_{1}\ldots  d_{k}).$ 
\end{itemize}
\end{itemize}
\end{teo1}

\section{Proof of Theorem 1}

Let $Z_{\ast}(X)$ denote the group of cycles on $X$ and
$$
[\,\,\,\,]: Z_{\ast}(X) \longrightarrow H_{\ast}(X;\CC)
$$
\noindent the cycle map. Since $X$ is compact, we observe that Borel--Moore homology coincides with singular homology. Considering $V$ as a $(2n-2k)$-cycle on $X$, we denote by $[V] \in  H_{2n-2k}(X;\CC)$ the
fundamental class determined by $V$.

By hypothesis, there is a holomorphic section of $N$, $s: X \longrightarrow N$, such that the zero set $Z(s)$ defines $V$. Since $s$ is a regular section, we have (see \cite[Proposition 14.1]{WF})
\begin{eqnarray}\label{xx0}
Z(s) = X\frown c_k(N)  \in Z_{2n-2k}(X),
\end{eqnarray}
\noindent considering $X$ as a $2n$-cycle in $Z_{2n}(X)$. By applying the cycle map in (\ref{xx0}), we obtain
\begin{eqnarray}\label{xx1}
[V] = [Z(s)] = [X \frown c_k(N) ] \in H_{2n-2k}(X;\CC).
\end{eqnarray}
\noindent On the other hand, for all $2n$-cycle $\alpha \in Z_{2n}(X)$, we know that (see, for example, \cite[Proposition 19.1.2]{WF})
$$
[\alpha \frown c_k(N) ] = [\alpha] \frown c_k(N) \in H_{2n-2k}(X;\CC),
$$ 
\noindent in particular, for $\alpha = X$, we obtain
$$
[X \frown c_k(N)] = [X] \frown c_k(N),
$$
\noindent and it follows from (\ref{xx1}) that
$$
[V] = [X] \frown c_k(N)  \in H_{2n-2k}(X;\CC).
$$
\noindent Hence, the $k$-th Chern class of $N$ is Poincaré dual to
the fundamental class of $V$. Thus, for a homogeneous symmetric polynomial
 $\varphi$ of degree $d$, we have the relation
\begin{eqnarray}\label{xx20}
[V] \frown \varphi(N) = [X] \frown \left(\varphi(N) \wedge c_k(N)\right).
\end{eqnarray}

According to \cite[Ch. VI, Theorem 6.9]{Suw2}, the sum of residues is given by 

$$
\sum_{\lambda=1}^r \mbox{Res}_{\varphi}(\fol,N,S_{\lambda}) = [V] \frown \varphi(N),
$$ 

\noindent hence, we get 
$$
\sum_{\lambda=1}^r \mbox{Res}_{\varphi}(\fol,N,S_{\lambda}) = [X] \frown \left(\varphi(N) \wedge c_k(N)\right).
$$ 

\noindent Thus, item (a) is proved.

If $V = \bigcap_{l =1}^kD_{l}$ is a complete intersection  of $k$ hypersurfaces $D_{l} \subset X$, then $N = \bigoplus^k_{\lambda=1}L(D_{\lambda})$, where $L(D_{\lambda})$ denotes the line bundle of the divisor $D_{\lambda}$. Also, $V$ coincides with the zero set $Z(s)$ of the section $s: X \longrightarrow N = \bigoplus^k_{\lambda=1}L(D_{\lambda})$ defined as follows: covering $X$ by open sets $U_{\alpha}$ so that each hypersurface $D_{\lambda}$ is defined by an equation $f_{\lambda ,\alpha} = 0$, in $U_{\alpha}$, we have that each holomorphic bundle $L(D_{\lambda})$ admits a natural holomorphic section determined by the collection $\{f_{\lambda, \alpha}\}$. The collection $\{\bigoplus^k_{\lambda=1}f_{\lambda \alpha}\}$ determines the holomorphic section $s$ of $\bigoplus^k_{\lambda=1}L(D_{\lambda})$ whose zero set $Z(s)$ is exactly $V = \bigcap_{l =1}^kD_{l}$.

For $\varphi = c_{n-k}$, we have
\begin{eqnarray}\nonumber
\varphi(N) \wedge c_k(N) &=& c_{n-k}(\bigoplus^k_{\lambda=1}L(D_{\lambda})) \wedge c_k(\bigoplus^k_{\lambda=1}L(D_{\lambda}))\\\nonumber
&=& \hspace{-0.3cm} \sum_{1\leq l_1 <\ldots <l_{n-k} \leq k}\hspace{-0.8cm} c_1(L(D_{l_1}))... c_1(L(D_{l_{n-k}}))\prod^k_{l=1}c_1(L(D_l))
\end{eqnarray}

\noindent and item (b)(i) is proved. 

Also, if $X = \PP^n$ and $deg(D_l) = d_l$, we get
\begin{eqnarray}\nonumber
c_1(L(D_l)) &=& c_1(\OO(d_l))\\\nonumber
&=& d_lc_1(\OO(1)),
\end{eqnarray}
\noindent and, therefore,

{\small{
 
\begin{eqnarray}\nonumber
&&[\PP^n] \frown \hspace{-0.3cm} \sum_{1\leq l_1 <\ldots <l_{n-k} \leq k}\hspace{-0.8cm} c_1(L(D_{l_1}))... c_1(L(D_{l_{n-k}}))\prod^k_{l=1}c_1(L(D_l))= \\\nonumber\\\nonumber && = [\PP^n] \frown \hspace{-0.3cm} \sum_{1\leq l_1 <\ldots <l_{n-k} \leq k}\hspace{-0.8cm} d_{l_1}... d_{l_{n-k}} c_1(\OO(1))^{n-k}d_{1}... d_{k}c_1(\OO(1))^{k}= \\\nonumber\\\nonumber && = 
 \hspace{-0.3cm} \sum_{1\leq l_1 <\ldots <l_{n-k} \leq k}\hspace{-0.8cm} (d_{l_1}... d_{l_l})(d_{1}... d_{l}) \left([\PP^n] \frown c_1(\OO(1))^{n}\right)= \\\nonumber\\\nonumber && = 
 \hspace{-0.3cm} \sum_{1\leq l_1 <\ldots <l_{n-k} \leq k}\hspace{-0.8cm} (d_{l_1}... d_{l_l})(d_{1}... d_{l}) 
\end{eqnarray}
}}
\noindent where in the last step we are using the relation
$$
[\PP^n] \frown c_1(\OO(1))^{n} = 1.
$$

\noindent Thus, item(c)(i) is proved.

On the other hand, for $\varphi = c_1^{n-k}$, we have
\begin{eqnarray}\nonumber
\varphi(N) \wedge c_k(N) &=& (c_1(\bigoplus^k_{\lambda=1}L(D_{\lambda})))^{n-k} \wedge c_k(\bigoplus^k_{\lambda=1}L(D_{\lambda}))\\\nonumber
&=& \hspace{-0.3cm} \left(\sum_{l=1}^{k} c_1(L(D_{l}))\right)^{n-k}\prod^k_{l=1}c_1(L(D_l))
\end{eqnarray}

\noindent and the item (b)(ii) is proved. If $X = \PP^n$ and $deg(D_l) = d_l$, we get

{\small{
 
\begin{eqnarray}\nonumber
&&[\PP^n] \frown \left(\sum_{l=1}^{k} c_1(L(D_{l}))\right)^{n-k}\prod^k_{l=1}c_1(L(D_l)) = \\\nonumber  \\\nonumber && = [\PP^n] \frown \left(\left(\sum_{l=1}^{k} d_l\right)c_1(\OO(1))\right)^{n-k}d_{1}... d_{k}c_1(\OO(1))^{k} = \\\nonumber  \\\nonumber && = \left(d_1 +\ldots + d_k\right)^{n-k}d_{1}\cdot\ldots \cdot d_{l} \left([\PP^n] \frown c_1(\OO(1))^{n}\right)= \\\nonumber  \\\nonumber && =  \left(d_1 +\ldots + d_k\right)^{n-k}d_{1}\cdot\ldots \cdot d_{l} 
\end{eqnarray}
}}

\noindent Thus, item(c)(ii) is proved.

\noindent $\square$

\section{Applications}
A. Lins Neto \cite{AL1} considered foliations by curves on $\PP^2$
that leave invariant one-dimensional irreducible subvarieties, and obtained a formula that relates the total sum of the Camacho-Sad index of the foliation with the degree and other characteristics of the invariant subvariety. More precisely: 
\begin{teo1234}[A. Lins Neto, \cite{AL1}]
Let $\fol$ be a one-dimensional foliation on $\PP^2$ and $D$ be a compact algebraic irreducible solution of $\fol$. Then
\begin{eqnarray}\label{ff00}
CS(\fol, D) = 3\,deg(D) - \chi(D) + \sum_{i=1}^r\mu_{x_i}(D),
\end{eqnarray}
\noindent where $x_1, x_2, \ldots, x_r$ are all singularities of $D$, $\mu_{x_i}(D)$ denotes the Milnor number of $D$ at $x_i$ and $CS(\fol, D)$ denotes the total sum of the Camacho-Sad index of $\fol$.
\end{teo1234}

One of the problems proposed by A. Lins Neto was to generalize formula (\ref{ff00}) for higher dimensions (see \cite{AL1}, {\it 1$^{st}$ problem}). As an application of Theorem~1, we obtain the following generalization of formula~(\ref{ff00}) for arbitrary dimension, which provides an answer to the proposed question:

\begin{teo2} \label{teo00}
Let $\fol$ be a one-dimensional foliation on $\PP^n$ and $D\subset \PP^n $ be a compact hypersurface such that $\fol$ leaves $D$ invariant. Then
{{\small
\begin{eqnarray}\label{ff01}
CS(\fol, D) =  \sum^{n-2}_{j=0}\binom{n+1}{n-1-j}(-1)^{n+j}deg(D)^{j+1}  - (-1)^n\chi(D) + \sum_{\lambda=1}^r\mu_{S_{\lambda}}(D),
\end{eqnarray}}}

\noindent where $S_1, S_2, \ldots, S_r$ are all connected components of $Sing(D)$, $\mu_{S_{\lambda}}(D)$ denotes the Milnor number of $D$ at ${S_{\lambda}}$ and $CS(\fol, D)$ denotes the total sum of residues of type Camacho-Sad of $\fol$.
\end{teo2}

\begin{remark}
\noindent{\bf (a)} A. Lins Neto \cite{AL1}  considered the invariant hypersurface $D$ with isolated singularities. In Theorem 2, the invariant hypersurface $D$ may have singularities of arbitrary type.
\\

\noindent{\bf (b)} To describe the total sum of the Camacho-Sad index, A. Lins Neto \cite{AL1} considered the first Chern polynomial $c_1$. Accordingly, in Theorem 2, we consider the Camacho-Sad index defined by the homogeneous symmetric polynomial $\varphi = c_1^{n-1}$.
\\

\noindent{\bf (c)} In the formula (\ref{ff01}), the sum $\displaystyle \sum^{n-2}_{j=0}\binom{n+1}{n-1-j}(-1)^{n+j}deg(D)^{j+1}$ reduces to $3\,deg(D)$, when $n =2$. Therefore, formula (\ref{ff01}) 
coincides with formula (\ref{ff00}), if $n=2$ and $D$ has only isolated singularities.
\\

\noindent{\bf (d)} D. Lehmann and T. Suwa \cite{LS00} also obtained a generalization of formula (\ref{ff00}).
\end{remark}

The proof of Theorem 2  will be presented in Subsection $\ref{s00}$. In addition to using Theorem 1, we first need to obtain a new version of the adjunction formula (see Proposition \ref{teo0}), which we will discuss next.

\subsection {On adjunction formulas}
Adjunction formulas are fundamental tools that arise naturally in the context of the classification theory of algebraic varieties. In particular, in the study of projective algebraic varieties, an effective strategy for exploring the classification spectrum is to consider the canonical class as a central invariant. The analysis of the canonical class captures essential birational properties of the variety. When comparing two varieties with similar structures, it is common to seek expressions that relate their canonical divisors. These expressions, known as adjunction formulas, allow one to describe the canonical divisor of a subvariety or a fiber in terms of the ambient canonical divisor and other associated geometric invariants.

If $C$ is a compact (possibly) singular curve in a complex compact surface $M$, a classical adjunction formula (due to K. Kodaira \cite{kod}) says that
\begin{eqnarray} \label{f0}
\chi(C) = - (K_M +C).C + \sum_{i=1}^r\mu_{x_i}(C)
\end{eqnarray}
\noindent where $x_1, x_2, \ldots, x_r$ are all singularities of $C$ and $K_{M}$ is the canonical divisor of $M$. Since $K_M.C = - [C]\frown c_1(TM)$, formula (\ref{f0}) can be written in the following terms
\begin{eqnarray} \label{f1}
\chi(C) = [C]\frown c_{n-1}(TM - L(C)) + \sum_{i=1}^r\mu_{x_i}(C).
\end{eqnarray}

More recently, in the work of J. Seade and T. Suwa  \cite{SeaSuw3}, the above formula has been generalized to the higher dimensional case:

\begin{teo1234} [Seade-Suwa,  \cite{SeaSuw3}] \label{teo000}
Let $X$be an n-dimensional complex compact manifold and $V\subset X$ a local complete intersection of codimension $k$ defined by a section of a holomorphic vector bundle $N$. If $V$ is compact and has only isolated singularities $x_1, x_2, \ldots, x_r$, then
\begin{eqnarray}\label{for1}
\chi(V) = [V]\frown c_{n-k}(TX - N) - (-1)^{n-k}\sum_{i =1}^r\mu_{x_{i}}(V).
\end{eqnarray}
\end{teo1234}

\begin{remark}
\noindent T. Suwa \cite{Suw3} proved a generalization of (\ref{f1}), where the complex surface $M$ is not assumed to be compact. 
\end{remark}

In formula (\ref{for1}), the variety $V$ has only isolated singularities. We obtain the following result, which provides an adjunction formula for a hypersurface $D$ that may have singularities of any type:

\begin{prop} \label{teo0}
Let $X$ be an n-dimensional complex compact manifold and $D\subset X$ a compact (possibly) singular hypersurface. Then
\begin{eqnarray}\label{for80}
\chi(D) = [D]\frown c_{n-k}(TX - L(D)) - (-1)^{n-1}\sum_{\lambda=1}^r\mu_{S_{\lambda}}(D)
\end{eqnarray}
\noindent where $S_1, S_2, \ldots, S_r$ are all connected components of $Sing(D)$, $\mu_{S_{\lambda}}(D)$ denotes the Milnor number of $D$ at ${S_{\lambda}}$ and $L(D)$ is the line bundle of $D$.
\end{prop}

\noindent {\bf Proof.} According to \cite[Proposition 1.6]{adam}, we have
$$
[X]\frown c_n(TX^{\vee} \otimes L(D)) = \sum_{\lambda =1}^r\mu_{S_{\lambda}}(D) - (-1)^n\chi(D)  + (-1)^n\chi(X).
$$

\noindent By using the definition of Chern class of tensor product, we can calculate the top
Chern class of $TX^{\vee} \otimes L(D)$ as follows:

{{\small
\begin{eqnarray} \nonumber
c_n(TX^{\vee} \otimes L(D)) &=& \sum_{j=0}^nc_j(TX^{\vee})c_1(L(D))^{n-j}\\\nonumber &=& c_n(TX^{\vee}) + \left(\sum_{j=0}^{n-1}c_j(TX^{\vee})c_1(L(D))^{n-1-j}\right)c_1(L(D))\\\nonumber &=& c_n(TX^{\vee}) + \left(\sum_{j=0}^{n-1}(-1)^jc_j(TX)c_1(L(D))^{n-1-j}\right)c_1(L(D))\\\nonumber &=& c_n(TX^{\vee}) + \left(\sum_{j=0}^{n-1}(-1)^{n-1}(-1)^{n-1-j}c_j(TX)c_1(L(D))^{n-1-j}\right)c_1(L(D))\\\nonumber &=& c_n(TX^{\vee}) + \left((-1)^{n-1}\sum_{j=0}^{n-1}c_j(TX)c_1(L(D)^{\vee})^{n-1-j}\right)c_1(L(D))\\\nonumber &=& c_n(TX^{\vee}) + (-1)^{n-1}c_{n-1}(TX - L(D))c_1(L(D))
\end{eqnarray}
}}

\noindent Thus, using the fact that $c_1(L(D))$ is the Poincaré dual of the fundamental class of $D$, we obtain
$$
[X]\frown c_n(TX^{\vee}) + (-1)^{n-1}[D]\frown c_{n-1}(TX - L(D)) = \sum_{\lambda =1}^r\mu_{S_{\lambda}}(D) - (-1)^n\chi(D) + (-1)^n\chi(X).
$$
\noindent and formula (\ref{for80}) follows from the classical Gauss Bonnet formula.

\noindent $\square$

\subsection {Proof of Theorem 2} \label{s00} 

Since 
$$
c_{n-1}(T\PP^n - L(D)) = \sum_{j=0}^{n-1}(-1)^jc_{n-1-j}(T\PP^n)c_1(L(D))^j
$$
\noindent we get
$$
c_1(L(D))^n = (-1)^{n+1}\left(c_{n-1}(T\PP^n - L(D))c_1(L(D))+ \sum^{n-2}_{j=0}(-1)^jc_{n-1-j}(T\PP^n)c_1(L(D))^{j+1}   \right).
$$
\noindent Thus, by item (b)(ii) of Theorem 1, we obtain

{\small{
\begin{eqnarray}\nonumber
 && \hspace*{-0.65 cm} CS(\fol,D) = [\PP^n] \frown c_1(L(D))^n\\\nonumber
&& \hspace*{-0.65 cm} =  (-1)^{n+1}\left([\PP^n] \frown c_{n-1}(T\PP^n - L(D))c_1(L(D)) + [\PP^n] \frown \sum^{n-2}_{j=0}(-1)^jc_{n-1-j}(T\PP^n)c_1(L(D))^{j+1} \right)\\\nonumber
&& \hspace*{-0.65 cm} =  (-1)^{n+1}\left([D] \frown c_{n-1}(T\PP^n - L(D)) + [\PP^n] \frown  \sum^{n-2}_{j=0}(-1)^jc_{n-1-j}(T\PP^n)c_1(L(D))^{j+1} \right)
\end{eqnarray}
}}

\noindent where in the last equality we have used the fact that $c_1(L(D))$ is Poincaré dual to the fundamental class of $D$. 

Now, since $c_1(L(D)) = deg\,(D)c_1(\OO_{\PP^n}(1))$ and $$c_{n-1-j}(\PP^n) = {n+1 \choose n-1-j}c_1(\OO_{\PP^n}(1))^{n-1-j},\,\,\, \mbox{for all $j=0,\ldots,n-2$,}$$
\noindent we have {\small{
\begin{eqnarray}\nonumber
CS(\fol, D) &=& \hspace*{-0.25 cm}  (-1)^{n+1}\left([D] \frown c_{n-1}(T\PP^n - L(D)) +  \sum^{n-2}_{j=0}\binom{n+1}{n-1-j}(-1)^{j+1}deg(D)^{j+1}\right)\\\nonumber &=& \hspace*{-0.25 cm}  (-1)^{n+1} \left([D] \frown c_{n-1}(T\PP^n - L(D))\right) +  \sum^{n-2}_{j=0}\binom{n+1}{n-1-j}(-1)^{n+j}deg(D)^{j+1}
\end{eqnarray}}}
\noindent and formula (\ref{ff01}) follows from adjunction formula (\ref{for80}).

\noindent $\square$

\begin{cor} \label{teo00009}
Let $\fol$ be a one-dimensional foliation on $\PP^n$ and $D\subset \PP^n$ be a hypersurface of degree $k$. If $\fol$ leaves $D$ invariant, then the numbers $\mu(D)=\sum_{\lambda=1}^r\mu_{S_{\lambda}}(D)$ and $\chi(D)$ satisfy the following compatibility condition
\begin{eqnarray}\label{ff000001}
\mu(D) - (-1)^n\chi(D) \equiv 0\, (mod\,\,\, k).
\end{eqnarray}
\end{cor}

\noindent {\bf Proof.}  It follows from item $(c)(ii)$ of Theorem 1 and formula (\ref{ff01}) of Theorem 2 that $k$ is a zero of the polynomial 
$$
P(X) = X^n+b_{n-1}X^{n-1}+\ldots +b_1X + b_0,
$$
\noindent where $\displaystyle b_j = - \binom{n+1}{n-j}(-1)^{n-1+j}$, for $j=1,\ldots,n-1$, and $b_0 = (-1)^n\chi(D) - \mu(D)$. Since $P(X)$ is a polynomial with integer coefficients, we have that $k$ is a factor of the constant term $b_0 = (-1)^n\chi(D) - \mu(D)$, hence we obtain (\ref{ff000001}).\\
\noindent $\square$

\begin{cor}\label{corvv}
Let $\fol$ be a one-dimensional foliation on $\PP^2$, and let 
$D \subset \PP^2$ be an irreducible algebraic solution of $\fol$. 
If $D$ has a singular point, then the degree of $D$ is at least three.
\end{cor}

\noindent {\bf Proof.} Suppose $D$ has a singular point, that is, $Sing(D) \neq \emptyset$. In particular, we have that $\mu(D) > 0$.  Assume $\deg(D) \in \{1,2\}$, which will be shown to lead to a contradiction. Of course, formula (\ref{ff01}) gives
$$
-2 =  - \chi(D) + \mu(D).
$$
\noindent Since $D$ is irreducible, we have that $\chi(D)\leq 2$ and, consequently, $\mu(D) = 0$, a contradiction.\\
\noindent $\square$

The following example shows that the hypothesis that $D$ is irreducible is an essential condition in Corollary \ref{corvv}.

\begin{exe}
Let $\PP^2$ be the complex projective space of dimension 2, with homogeneous coordinates $(X_0:X_1:X_2)$, and let $\fol$ be the holomorphic foliation on $\PP^2$ defined by the vector field (affine coordinate $X_0 \neq 0$)
$$
v =  x_1\displaystyle\frac{\partial}{\partial x_1} + x_2\displaystyle\frac{\partial}{\partial x_2}.
$$

\noindent The hypersurface $D \subset \mathbb{P}^2$, defined by $X_1^2 + X_2^2 = 0$,
\noindent is a reducible curve of degree less than 3 that is $\fol$-invariant, whose singular set is $Sing(D) = \{(1:0:0)\}$. 

In the affine chart $X_0 \neq 0$, with coordinates $(x_1,x_2) = (\frac{X_1}{X_0}, \frac{X_2}{X_0}),$
the local defining function of $D$ is $f(x_1,x_2) = x_1^2 + x_2^2$. Since $f$ a {\it Pham-Brieskorn} polynomial of degree $deg(D) = 2$, we have that $\mu_{(1:0:0)}(D) = 1$. Also, since $D$ is the union of two projective lines meeting at 
the singular point $(1:0:0)$, we have that $\chi(D) = 3$.

Hence, 
$$
3\,deg(D) - \chi(D) + \sum_{i=1}^r\mu_{x_i}(D) = 4 = (deg(D))^2
$$

\noindent according to formula (\ref{ff01}).

\end{exe}

\begin{cor}\label{corrr}
Let $\fol$ be a one-dimensional foliation on $\PP^2$, and let 
$D \subset \PP^2$ be an irreducible algebraic solution of $\fol$. Then
\begin{eqnarray}
\# \Sing(D) \leq deg(D)(\deg(\fol) - 1) +2.
\end{eqnarray}

\noindent In particular, 
\begin{eqnarray}\label{desss}
deg(D) \leq deg(\fol) +2,\,\,\,\mbox{if $D$ has at most nodal singularities,}
\end{eqnarray}
\noindent and
\begin{eqnarray}
\,\,\,\,\,\,\,\,\,\,\,\# \Sing(D) \leq deg(\fol)(\deg(\fol) + 1), \mbox{if $D$ is non-dicritical with respect to $\fol$}.
\end{eqnarray}
\end{cor}

\begin{remark}
Inequality  (\ref{desss}) recovers a result due to D. Cerveau and A. Lins Neto \cite[Ch. 3, Theorem 1]{Alc}
\end{remark}

\noindent {\bf Proof of Corollary \ref{corrr}:} Let $\fol$ be a one-dimensional foliation on $\PP^2$ of degree $deg(\fol) = d$, and let 
$D \subset \PP^2$ be an irreducible algebraic solution of $\fol$ of degree $deg(D) = k$.  It follows from \cite[Theorem C]{DSM} and \cite[Theorem 2]{DSM2} that
\begin{eqnarray}\label{foooo}
\# \Sing(D) \leq \sum Sch_p(\fol,D) = 2k -k^2 +dk + \sum_{p\in Sing(D)}\mu_p(D)
\end{eqnarray}
\noindent where $\sum Sch_p(\fol,D)$ denotes the sum of the local Schwarz-index of $\fol$. On the other hand, since $D$ irreducible, formula (\ref{ff01}) tell us that 
$$
\sum_{p\in Sing(D)}\mu_p(D)  = \chi(D) -3k +k^2 \leq 2 -3k +k^2
$$
\noindent and, consequently, 
\begin{eqnarray}\label{dfg}
\# \Sing(D) \leq k(d-1) +2.
\end{eqnarray}
Now, if $D$ has at most nodal singularities, then each singular point $p\in Sing(D)$ has
multiplicity $1$. Thus, it follows from (\ref{foooo}) that
$$
0 \leq k(2+d - k),
$$
\noindent or, equivalently,
$$
k \leq d+2.
$$

Finally, if $D$ is non-dicritical with respect to $\fol$, then according to \cite[Proposition 6]{Bru} we have $k \leq d+2$. By using this inequality in (\ref{dfg}), we obtain
$$
k\leq d(d+1).
$$

\noindent $\square$

\subsection {Formulas for Milnor number} \label{s0xc} 

Formulas $(\ref{ff01})$ and $(\ref{for80})$ involve the Milnor number of a hypersurface with nonisolated singularities, a concept introduced by A. Parusinski in \cite{adam}. In this section, under mild conditions on the singular locus of the hypersurface, we provide  explicit formulas for this invariant (Theorem \ref{parv}, Corollary \ref{parvv}).

In this section, let $D \subset X$ be a hypersurface with non-isolated singularities on an $n$-dimensional compact complex manifold $X$. Let $S$ be a connected component of $\Sing(D)$, assumed to be a smooth local complete intersection of codimension $k$, given as the zero set of a section of a holomorphic vector bundle $N$ on $X$.

\begin{teo}\label{parv} If the pair $(D-Sing(D), S)$ satisfies Whitney's conditions and $\mu$ denotes $(n-k)$-th sectional Milnor number of $D$, then 
{\small {
\begin{eqnarray}\nonumber
 \mu_{S}(D) = \mu\,.\int_{X}\hspace{-0.04 cm}\sum_{t=0}^{n-k} \hspace{-0.0 cm} \sum_{j=0}^{n-k-t}\hspace{-0.0 cm}\sum_{i=1}^j\sum_{\mid L\mid = j}\hspace{-0.09 cm}(-1)^{n-k-t+i}c_{n-k-t-j}(TX)c_{L}(N)c_{k}(N) c_1(L(D))^t
\end{eqnarray}
}}

\noindent where $L:=L^i = (l_1,\ldots,l_i)\in (\NN^{\ast})^i$, with  $l_1+\ldots +l_i = j$, is an $i$-dimensional multi-index and $\displaystyle c_{L}\left(N\right) := c_{l_1}\left(N\right)\cdots c_{l_i}\left(N\right)$.
\end{teo}

\begin{remark}
\noindent{\bf (a)} It is knoww that $(D-Sing(D), S)$ satisfies Whitney's conditions iff for each $k=1,\ldots, n-k$, the $k$-th sectional Milnor  number of $D$ is constant on $S$ (see \cite{TT2}). In such a situation, $\mu$ is defined and can be calculated as the $(n-k)$-th sectional Milnor at any point of $S$. 
\\

\noindent{\bf (b)} The case where $S$ is a local complete intersection defined
as the zero set of a section of a vector bundle is the case, for example, where $S$ is a global complete intersections or a complete intersections in projectives
spaces.
\\
\end{remark}

\noindent {\bf Proof of Theorem \ref{parv}:} 
By \cite[Proposition 1.7]{adam}, we have 
$$
\mu_S(D) = \mu\,.\int_Sc_{n-k}(TS^{\vee} \otimes L(D)).
$$

\noindent The $(n-k)$-th Chern class of tensor product $TS^{\vee} \otimes L(D)$ can be calculated as follows:
\begin{eqnarray}\nonumber
\mu_S(D) &=& \mu\,.\int_Sc_{n-k}(TS^{\vee} \otimes L(D))\\\nonumber
&=& \mu\,.\int_S\sum_{t=0}^{n-k}(-1)^{n-k-t}c_{n-k-t}(TS)c_1(L(D))^t\\\nonumber
&=& \mu\,.\int_S\sum_{t=0}^{n-k}(-1)^{n-k-t}c_{n-k-t}(TX - N)c_1(L(D))^t
\end{eqnarray}
\noindent where in the last step, we make use of the fact that $TS = (TX - N)|_S$, as the 
component $S$ is a local complete intersection defined as
the zero set of a section of $N$. On the other hand, using the fact that $c_{k}(N)$ is the Poincaré dual of the fundamental class of $S$, we obtain
$$
\mu_S(D) = \mu\,.\int_X\sum_{t=0}^{n-k}(-1)^{n-k-t}c_{n-k-t}(TX - N)c_1(L(D))^tc_k(N).
$$
\noindent Now, since $c_{\ast}(TX - N) = c_{\ast}(TX)c_{\ast}^{-1}(N)$, we get 
$$
c_{n-k-t}(TX - N) = \sum_{j=0}^{n-k-t}\hspace{-0.0 cm}\sum_{i=1}^j\sum_{\mid L\mid = j}\hspace{-0.0 cm}(-1)^{i}c_{n-k-t-j}(TX)c_L(N)
$$
\noindent where $L:=L^i = (l_1,\ldots,l_i)\in (\NN^{\ast})^i$, with  $l_1+\ldots +l_i = j$ and $\displaystyle c_{L}\left(N\right) := c_{l_1}\left(N\right)\cdots c_{l_i}\left(N\right)$, and, consequently,

{\small {
\begin{eqnarray}\nonumber
 \mu_{S}(D) = \mu\,.\int_{X}\hspace{-0.04 cm}\sum_{t=0}^{n-k} \hspace{-0.0 cm} \sum_{j=0}^{n-k-t}\hspace{-0.0 cm}\sum_{i=1}^j\sum_{\mid L\mid = j}\hspace{-0.09 cm}(-1)^{n-k-t+i}c_{n-k-t-j}(TX)c_{L}(N)c_{k}(N) c_1(L(D))^t.
\end{eqnarray}
}}

\noindent $\square$

\begin{cor}\label{parvv}
Setting as in Theorem \ref{parv}. Suppose $X = \PP^n$ and $S =  \bigcap_{\lambda =1}^{k}D_{\lambda}$ a complete intersection of $k$ hypersurfaces, then
{\small{
\begin{eqnarray}\nonumber
 \mu_{S}(D) = \mu\,.\prod_{\lambda = 1}^k \hspace{-0.1cm} deg(D_{\lambda})\hspace{-0.04 cm}\sum_{t=0}^{n-k}\hspace{-0.0 cm} \sum_{j=0}^{n-k-t}\hspace{-0.0 cm}\sum_{i=1}^j\sum_{\mid L\mid = j}\hspace{-0.0 cm}(-1)^{n-k-t+i}\binom{n+1}{n-k-t-j}deg_L(S)deg(D)^t\hspace{-0.1 cm}
\end{eqnarray}
}}

\noindent where $L:=L^i = (l_1,\ldots,l_i)\in (\NN^{\ast})^i$, with  $\mid L\mid = l_1+\ldots +l_i$ and 
{\small{
$$\displaystyle deg_L(S) := \hspace{-0.12 cm} \left(\displaystyle\sum_{1\leq e_{1}^1<\ldots<e_{l_1}^1 \leq k}\hspace{-0.7 cm}deg(D_{e_{1}^1}) \cdot\ldots\cdot deg(D_{e_{l_1}^1})\hspace{-0.05 cm}\right)\cdot \ldots \cdot \left(\displaystyle\sum_{1\leq e_{1}^i<\ldots<e_{l_i}^i \leq k}\hspace{-0.7 cm}deg(D_{e_{1}^i}) \cdot \ldots \cdot deg(D_{e_{l_i}^i})\right).$$
}}
\end{cor}

\noindent {\bf Proof.} Let  $S =  \bigcap_{\lambda =1}^{k}D_{\lambda}$ be is a complete intersection of $k$ hypersurfaces. In this case, $S$ is defined as
the zero set of a section of the holomorphic vector $N = \bigoplus_{\lambda =1}^{k} L(D_{\lambda})$, and, consequently, we have that
$$
c_k(N) = \prod_{\lambda = 1}^k \hspace{-0.1cm} deg(D_{\lambda})\OO_{\PP^n}(1)^k
$$
\noindent and
{\footnotesize  {$$
c_L(N) = \hspace{-0.12 cm} \left(\displaystyle\sum_{1\leq e_{1}^1<\ldots<e_{l_1}^1 \leq k}\hspace{-0.7 cm}deg(D_{e_{1}^1}) \cdot\ldots\cdot deg(D_{e_{l_1}^1})\hspace{-0.05 cm}\right)\cdot \ldots \cdot \left(\displaystyle\sum_{1\leq e_{1}^i<\ldots<e_{l_i}^i \leq k}\hspace{-0.7 cm}deg(D_{e_{1}^i}) \cdot \ldots \cdot deg(D_{e_{l_i}^i})\right)\OO_{\PP^n}(1)^j.
$$}}

\noindent for each $L=L^i = (l_1,\ldots,l_i)\in (\NN^{\ast})^i$, with  $l_1+\ldots +l_i = j$. Thus, using formula of  the Theorem \ref{parv}, we obtain
{\footnotesize {
\begin{eqnarray}\nonumber
 \mu_{S}(D)\hspace{-0.3cm} &=& \hspace{-0.3cm} \mu\,.\int_{X}\hspace{-0.04 cm}\sum_{t=0}^{n-k} \hspace{-0.0 cm} \sum_{j=0}^{n-k-t}\hspace{-0.0 cm}\sum_{i=1}^j\sum_{\mid L\mid = j}\hspace{-0.09 cm}(-1)^{n-k-t+i}c_{n-k-t-j}(TX)deg_L(S)\prod_{\lambda = 1}^k \hspace{-0.1cm} deg(D_{\lambda}) c_1(L(D))^t\OO_{\PP^n}(1)^{j+k}\\\nonumber
&=& \hspace{-0.3cm} \mu\,.\prod_{\lambda = 1}^k \hspace{-0.1cm} deg(D_{\lambda})\hspace{-0.04 cm}\sum_{t=0}^{n-k}\hspace{-0.0 cm} \sum_{j=0}^{n-k-t}\hspace{-0.0 cm}\sum_{i=1}^j\sum_{\mid L\mid = j}\hspace{-0.0 cm}(-1)^{n-k-t+i}\binom{n+1}{n-k-t-j}deg_L(S)deg(D)^t
\end{eqnarray}
}}

\noindent where in the last step we are using the relations 
$$
\displaystyle c_{n-k-t-j}(\PP^n) = \binom{n+1}{n-k-t-j}\OO_{\PP^n}(1)^{n-k-t-j}
$$ 

\noindent and 

$$
c_1(L(D))=deg(D)\OO_{\PP^n}(1).
$$

\noindent $\square$

\begin{remark}
\noindent{\bf (a)} Expression 
$$
\int_{S}c_{n-k}(TS) = \prod_{\lambda = 1}^k \hspace{-0.1cm} deg(D_{\lambda}) \sum_{j=0}^{n-k}\sum_{i=1}^j\sum_{\mid L\mid = j}(-1)^i\binom{n+1}{n-k-j}deg_L(S)
$$
\noindent where $L:=L^i = (l_1,\ldots,l_i)\in (\NN^{\ast})^i$, with  $\mid L\mid = l_1+\ldots +l_i$ and 
{\small{
$$\displaystyle deg_L(S) := \hspace{-0.12 cm} \left(\displaystyle\sum_{1\leq e_{1}^1<\ldots<e_{l_1}^1 \leq k}\hspace{-0.7 cm}deg(D_{e_{1}^1}) \cdot\ldots\cdot deg(D_{e_{l_1}^1})\hspace{-0.05 cm}\right)\cdot \ldots \cdot \left(\displaystyle\sum_{1\leq e_{1}^i<\ldots<e_{l_i}^i \leq k}\hspace{-0.7 cm}deg(D_{e_{1}^i}) \cdot \ldots \cdot deg(D_{e_{l_i}^i})\right).$$
}}

\noindent is the Euler characteristic of a (nonsingular) complete intersections $S =  \bigcap_{\lambda =1}^{k}D_{\lambda}$  of $k$ hypersurfaces in $\PP^n$.
\\

\noindent{\bf (b)} If $S =  \bigcap_{\lambda =1}^{n-1}D_{\lambda}$ is a complete intersections of $n-1$ hypersurfaces,   Corollary \ref{parvv} tells us that
\begin{eqnarray}\nonumber
\mu_{S}(D) = \mu\,.\prod_{\lambda = 1}^{n-1} \hspace{-0.1cm} deg(D_{\lambda})\left\{\sum_{\lambda =1}^{n-1}deg(D_{\lambda}) + deg(D) - (n+1)\right\}
\end{eqnarray}
\noindent or, equivalently,
$$
\mu_{S}(D) = \mu\,.\left\{\prod_{\lambda = 1}^{n-1} \hspace{-0.1cm} deg(D_{\lambda})deg(D) - \chi(S)\right\}.
$$
\end{remark}

\begin{exe}
Let $\PP^4 = \{(X_0:X_1:X_2:X_3:X_4)\}$ and $D$ the hypersurface in $\PP^4$ defined by
$$
X_0^2+ X_1^2 = 0.
$$
\noindent Then the singular set of $D$ is the (smooth) complete intersection of two hypersurfaces 
$$
Sing(D) = \{X_{0} = 0\} \cap \{X_{1} = 0\}
$$ 
\noindent In particular, $Sing(D)$ is given as the zero set of a section of the holomorphic vector bundle $N = L(D_0)\oplus L(D_1)$, where $D_0 = \{X_0 = 0\}$ and $D_1 = \{X_1 = 0\}$ .

We can apply the formula of Corollary \ref{parvv} to compute the Milnor number of $D$ at $S = \{(0:0:X_2:X_3:X_4)\}$: 
{\footnotesize{
\begin{eqnarray}\nonumber
\mu_{S}(D) = \mu  \{(-1)^2\displaystyle \binom{5}{2} + (-1)^3\binom{5}{1}\left[deg(D_0) + deg(D_1)\right]+(-1)^3\binom{5}{0}\left[deg(D_0) + deg(D_1)\right] +\\\nonumber (-1)^4\binom{5}{0}\left[deg(D_0) + deg(D_1)\right]^2+ (-1)^1\binom{5}{1}deg(D) + (-1)^2\binom{5}{0}deg(D)\left[deg(D_0) + deg(D_1) \right]+ \\\nonumber  (-1)^0\binom{5}{0}deg(D)^2 \}
\end{eqnarray}
}}

\noindent since $deg(D) = 2$ and $deg(D_0)= deg(D_1) =1$, we obtain $\mu_{S}(D) = 0$.
\end{exe}

\end{document}